\documentclass[11pt, reqno]{amsart}
\usepackage{amscd,amssymb}
\usepackage{amsfonts}
\usepackage{amsmath}
\usepackage{mathrsfs}
\usepackage{amscd} 
\usepackage{cite, tocvsec2}
\usepackage{amstext, amsthm}
\usepackage{shuffle}
\usepackage{amsopn,xspace,enumerate,multicol}
\usepackage{array}
\usepackage{mathtools}
\usepackage{cellspace}

\newcommand\lie[1]{\mathfrak{#1}}

\newcommand{\N}{\mathbb{N}}

\newcommand{\Q}{\mathbb{Q}}
\newcommand{\K}{\mathbb{K}}

\newtheorem{thm}{Theorem}
\newtheorem{cor}[thm]{Corollary}
\newtheorem{prop}[thm]{Proposition} 
 
\newtheorem{defn}[thm]{Definition} 

\newtheorem{exm}[thm]{Example}

\newtheorem*{rem}{Remark}

\newtheorem{lem}[thm]{Lemma}
\numberwithin{thm}{subsection}
\numberwithin{equation}{section}

\newcommand{\on}{\operatorname}
\newcommand{\ad}{\on{ad}} 
\newcommand{\ddiv}{\on{div}}
 
\newcommand{\qdiv}{\on{qdiv}}  
\newcommand{\tr}{\on{tr}}  
  
\newcommand{\qtr}{\on{qtr}}  

\newcommand{\grt}{\lie{grt}} 
\newcommand{\krv}{\lie{krv}} 
\newcommand{\tder}{\lie{tder}} 

\begin{document}

\title{Divergence and $q$-divergence in depth 2}
\author{Anton Alekseev, Anna Lachowska, Elise Raphael }

\date{\today} 

\begin{abstract}
The Kashiwara-Vergne Lie algebra $\mathfrak{krv}$ encodes symmetries of the Kashiwara-Vergne problem on the properties of the Campbell-Hausdorff series. It is conjectured that $\krv \cong \mathbb{K}t \oplus \mathfrak{grt}_1$, where $t$ is a generator of degree 1 and $\mathfrak{grt}_1$ is the Grothendieck-Teichm\"uller Lie algebra. In the paper, we prove this conjecture in depth 2. The main tools in the proof are the divergence cocycle and the representation theory of the dihedral group $D_{12}$. Our calculation is similar to the calculation by Zagier of the graded dimensions of the double shuffle Lie algebra in depth 2.

In analogy to the divergence cocycle, we define the super-divergence and $q$-divergence cocycles (here $q^l=1$) on Lie subalgebras of $\mathfrak{grt}_1$ which consist of elements with weight divisible by $l$. We show that in depth $2$ these cocycles have no kernel. This result is in sharp contrast with the fact that the divergence cocycle vanishes on
$[\mathfrak{grt}_1, \mathfrak{grt}_1]$.

\end{abstract}

\maketitle

\section{Introduction} 
The Kashiwara-Vergne problem in Lie theory \cite{KV} states the existence of an automorphism $F$ of the free Lie algebra with two generators $x$ and $y$ which satisfies the following properties:
\begin{itemize}
\item
The action of $F$ preserves the conjugacy classes of generators $x$ and $y$. That is, there exist $g(x,y)$ and $h(x,y)$ so as $F(x)=e^g xe^{-g}$ and $F(y)=e^h y e^{-h}$.

\item
The action of $F$ maps $x+y$ to the Campbell-Hausdorff series, $F(x+y)=\log(e^xe^y)$.

\item
The Jacobian cocycle of $F$  (see \cite{AT} for the precise definition) is of the form
$j(F)=f(\log(e^xe^y)) -f(x)-f(y)$.
\end{itemize}
The positive solution of the Kashiwara-Vergne problem (see \cite{AM} and \cite{AT}) implies the Duflo Theorem on the isomorphism of the center of the universal enveloping algebra and the ring of invariant polynomials, and various generalizations of algebraic and analytic aspects of the Duflo theorem. 

The Kashiwara-Vergne Lie algebra $\mathfrak{krv}$ encodes the symmetry of  the Kashiwara-Vergne problem. The corresponding pro-unipotent group acts freely and transitively on the space of solutions of the Kashiwara-Vergne problem. It has been shown in \cite{AT} that the Grothendieck-Teichm\"uller Lie algebra $\mathfrak{grt}_1$ injects into $\mathfrak{krv}$. It is conjectured that 
$$
\mathfrak{krv} \cong \mathbb{K} t \oplus \mathfrak{grt}_1,
$$
where $t$ is a generator of degree 1. Both Lie algebras $\mathfrak{krv}$ and $\mathfrak{grt}_1$ are graded by weight (the total number of $x$'s and $y$'s) and filtered by depth (the total number of $y$'s). 

We prove this conjecture in depth 2. That is, we show that $\mathfrak{krv}$ coincides with $\mathbb{K}t \oplus \mathfrak{grt}_1$ up to elements of depth 3 for arbitrary weight.
The main tools in the proof are the divergence cocycle ${\rm div}$ (the infinitesimal version of the Jacobian cocycle $j$) and the representation theory of the dihedral group $D_{12}$ acting on polynomials in two variables.
It is interesting to note that polynomials satisfying the same symmetry appeared previously in the work of 
Zagier (\cite{Z}, see also  
\cite{IhO}), in calculation of the graded dimensions of the double shuffle Lie algebra. 
L.~Schneps explained to us that her work \cite{S} implies that the double shuffle Lie algebra 
injects into $\krv$. This in combination with Zagier's result \cite{Z} could provide an alternative 
proof of the surjection of $\krv$ onto $\grt_1$ in depth 2. 

In analogy to the divergence cocycle, we introduce the super-divergence and $q$-divergence cocycles for $q$ a primitive root of unity of order $l$, $q^l=1$. They are defined on the Lie subalgebras of $\mathfrak{grt}_1$ which consist of elements with weight divisible by $l$.  The super-divergence cocycle plays an important role in Rouvi\`ere's theory of symmetric spaces (see \cite{R}).

We prove that in depth 2 the super-divergence and $q$-divergence cocycles have no kernel. This result is in sharp contrast with the fact that the divergence cocycle vanishes on $[\mathfrak{grt}_1, \mathfrak{grt}_1]$. We conjecture that super-divergence and $q$-divergence cocycles have no kernel in arbitrary depth.

\vskip 0.2cm

{\bf Acknowledgements.} We are grateful to B. Enriquez for very useful discussions and for his help with solving the functional equations for super-divergence and $q$-divergence. We are indebted to L. Schneps for explaining to us the relation of our results to the ones in her work \cite{S} and in the work of D. Zagier \cite{Z}. We would like to thank H. Furusho for bringing to our attention the results of Ihara and Takao. E.R. thanks P. de la Harpe for his generous help with questions of representation theory.

The research of A.A. was supported in part of the grant MODFLAT of the European Research Council (ERC), by the grants 152812 and 141329 of the Swiss National Science Foundation, and by the NCCR SwissMAP of the Swiss National Science Foundation. The research of E.R. was supported in part by the grant 152812 of the Swiss National Science Foundation.

\section{Cocycles on Lie algebras of tangential derivations.} 

\subsection{The Lie algebra of tangential derivations and the divergence cocycle.} 
Let $\lie{lie}_n = \lie{lie} (x_1, \ldots, x_n)$ denote the degree completion of the free Lie algebra over a field $\K$  on $n$ variables, and 
let ${\rm Ass}_n = U(\lie{lie}_n)$ be its universal enveloping algebra. Then ${\rm Ass}_n$ is a degree completion of the free associative algebra on 
$n$ generators.  We have 
\[ \lie{lie}_n = \prod_{k=1}^\infty \lie{lie}^k (x_1, \ldots x_n) , \] 
where $\lie{lie}^k(x_1, \ldots x_n)$ is spanned by the Lie words of length $k$. We will call $\lie{lie}^k (x_1, \ldots x_n)$ 
the graded component of $\lie{lie}_n$ of {\em weight} $k$. 

When $n=2$, we will use the variables $(x,y)$ instead of $(x_1, x_2)$. 

In \cite{AT}, the following graded vector space was defined: 
\[ \lie{tr}_n = {\rm Ass}_n^+ /\langle (ab-ba), \;\;\;  a,b \in {\rm Ass}_n \rangle ,\] 
where ${\rm Ass}_n^+$ is the augmentation ideal of ${\rm Ass}_n$. 

Recall the Lie algebra of tangential derivations on $\lie{lie}_n$: a derivation $u$ on $\lie{lie}_n$ is tangential 
if  there exist $a_1, \ldots, a_n \in \lie{lie}_n$ such that $u(x_i)= [x_i, a_i]$ for all $i = 1, \ldots n$. The action of $u$ on the generators 
completely determines the derivation. 

 The tangential derivations $\lie{tder}_n$ on $\lie{lie}_n$ form a Lie algebra. Namely, if 
$u = (a_1, \ldots a_n), \; v = (b_1, \ldots b_n) \in \lie{tder}_n$, then  $[u,v] = (c_1, \ldots c_n)$, where $c_k = u(b_k) - v(a_k) + [a_k, b_k]$ 
for all $k=1, \ldots n$.  

Every element $a \in {\rm Ass}_n$ has a unique decomposition 
\[  a =  a_0 + \sum_{k=1}^n (\partial_k a) x_k , \] 
where $a_0 \in \K$, and $(\partial_k a) \in {\rm Ass}_n$. 

The {\it divergence} map ${\rm div}: \lie{tder}_n  \longrightarrow \lie{tr}_n$ was defined in \cite{AT} as follows. Let $u = (a_1, \ldots a_n) 
\in \lie{tder}_n$. Then 
\[ {\rm div} (u)  = \sum_{k=1}^n {\rm tr} (x_k (\partial_k a_k))  \;\; \subset \;  \lie{tr}_n . \] 
The divergence is a cocycle on $\lie{tder}_n$, in the sense that $u \cdot {\rm div}(v) - v \cdot {\rm div}(u) = {\rm div}([u,v])$ for all 
$u, v \in \lie{tder}_n$. (see \cite{AT}, Proposition 3.6).

\subsection{Even tangential derivations and the superdivergence.} 
Similarly, we can define an even version of the above mentioned objects. 
Consider the subspace 
\[ \lie{lie}_n^{\rm {even}}  = \prod_{k \in 2 \cdot \N}  \lie{lie}^k(x_1, \ldots x_n) . \]  
This is a Lie subalgebra of $\lie{lie}_n$. Its universal enveloping algebra 
${\rm Ass}_n^{\rm {even}} = U(\lie{lie}^n)$ is spanned by the words of even length. 

 For ${\rm Ass}_n^{\rm{even}}$, the {\it supertrace}  graded vector space  can be defined as follows. 
\begin{defn}
The graded vector space $\lie{str}_n$ (the {\it supertrace}) is the quotient  
\[ \lie{str}_n = {\rm Ass}_n^{\rm{even} \;+} / \langle (x_i a + a x_i), \;\;\;\;\;\;  x_i \in \{ x_1, \ldots x_n\}, \; a \in {\rm Ass}_n^+  \rangle. \]     
\end{defn} 

Let ${\rm str}_n : {\rm Ass}_n^{\rm{even}} \longrightarrow \lie{str}_n$ denote the natural projection. 

\begin{exm} {\rm The graded component $\lie{str}_2^4 \subset \lie{str}_2$ is spanned by  \\ 
$\{  {\rm str}(x^4), \;  {\rm str}(x^3 y) = - {\rm str}(y x^3) = - {\rm str}(x^2 y x) = {\rm str}(xy x^2), \;  {\rm str}(x^2 y^2) = -{\rm str}(xy^2x) = {\rm str}(y^2x^2) = -{\rm str}(y x^2y), \; 
{\rm str}(xyxy) =-{\rm str}(yxyx),  \; {\rm str}(x y^3) = - {\rm str}(y^3 x) = - {\rm str}(y x y^2) = {\rm str}(y^2 x y),  \; {\rm str}(y^4) \}$. } 
\end{exm}      
 
Then we have a natural notion of even tangential derivations. 

 \begin{defn} The subspace $\lie{tder}_n^{\rm even} \subset \lie{tder}_n$ is spanned by the tangential derivations 
 $u=(a_1, a_2, \ldots a_n)$ such that $a_i \in \lie{lie}_n^{\rm even}$ for all $i = 1, \ldots n$. 
  \end{defn} 
\begin{prop} The even tangential derivations $\lie{tder}_n^{\rm even}$ form a Lie subalgebra of $\lie{tder}_n$. 
\end{prop} 
{\it Proof:}  If $a_k, b_k  \in \lie{tder}_n^{\rm even}$ for all $k = 1, \ldots n$,  then $c_k = u(b_k) - v(a_k) + [a_k, b_k] \in \lie{tder}_n^{\rm even}$ for all $k= 1, \ldots n$.  \qed 

The even tangential derivations preserve $\lie{lie}_n^{\rm even} \subset \lie{lie}_n$ and ${\rm Ass}_n^{\rm even} \subset {\rm Ass}_n$.  The action descends to $\lie{str}_n$ along the natural projection. If $u \in \lie{tder}_n^{\rm even}$, and $a \in  \lie{str}_n$, we will 
denote the result of this action by $u \cdot a \in \lie{str}_n$. 

Now we can define the superdivergence map for $\lie{tder}_n^{\rm even}$. 
\begin{defn} The {\em superdivergence} is the map ${\rm sdiv}: \lie{tder}_n^{\rm even}  \longrightarrow \lie{str}_n$ given by the formula 
\[ {\rm sdiv}(u = (a_1, \ldots a_n)) = \sum_{k=1}^n {\rm str} (x_k (\partial_k a_k))  \;\; \subset \;  \lie{str}_n . \] 
\end{defn} 

\begin{prop} The superdivergence is a cocycle on $\lie{tder}_n^{\rm even}$, in the sense that $u \cdot {\rm sdiv}(v) - v \cdot {\rm sdiv}(u) = {\rm sdiv}([u,v])$ for all $u, v \in \lie{tder}_n^{\rm even}$. 
\end{prop}  
{\it Proof:}  See Proposition \ref{qdiv} below. 
\qed

\subsection{The $q$-tangential derivations and the $q$-divergence.} 
More generally, let $q$ be a primitive $l$-th root of unity, $l \geq 1$, and assume that the field $\K$ contains $q$. Following the same logic as before, a notion of $q$-divergence can be defined. 

Let $\lie{lie}_n^{l \N}$ denote the Lie subalgebra of the free Lie algebra $\lie{lie}_n$ spanned 
by the Lie words of length divisible by $l$: 
\[ \lie{lie}_n^{l\N} = \prod_{k \in l \cdot \N} \lie{lie}^k (x_1, \ldots x_n) , \] 
and let ${\rm Ass}_n^{l\N}$ denote its universal enveloping algebra. 

\begin{defn} The graded vector space $\lie{qtr}_n$ (the $q$-trace) is the quotient 
\[ \lie{qtr}_n = {\rm Ass}_n^{l\N \;+} / \langle (x_i a -q a x_i), \;\;\;\;\;\;  x_i \in \{ x_1, \ldots x_n\}, \; a \in {\rm Ass}_n^+  \rangle. \]     
\end{defn} 

Let ${\rm qtr}_n : {\rm Ass}_n^{l\N} \longrightarrow \lie{qtr}_n$ denote the natural projection. 
\begin{exm} Let $l=3$. The graded component $\lie{qtr}_2^3 \subset \lie{qtr}_2$ is spanned by  \\ 
$\{ {\rm qtr}(x^3), \; {\rm qtr}(x^2 y) = q {\rm qtr}(x y x) = q^2 {\rm qtr}(y x^2), \; {\rm qtr}(x y^2) = q {\rm qtr}(y^2 x) = 
q^2 {\rm qtr}(y x y), 
\; {\rm qtr}(y^3) \}. $
\end{exm} 

 \begin{defn} The subspace $\lie{tder}_n^{l\N} \subset \lie{tder}_n$ is spanned by the tangential derivations 
 $u=(a_1, a_2, \ldots a_n)$ such that $a_i \in \lie{lie}_n^{l\N}$ for all $i = 1, \ldots n$. 
  \end{defn} 
  
The $q$-tangential derivations $\lie{tder}_n^{l\N}$ form a Lie subalgebra of $\lie{tder}_n$. The action of $\lie{tder}_n^{l\N}$ 
preserves $\lie{lie}_n^{l\N} \subset \lie{lie}_n$ and ${\rm Ass}_n^{l\N} \subset {\rm Ass}_n$, and descends to 
$\lie{qtr}_n$ along the natural projection.  If $u \in \lie{tder}_n^{l\N}$, and $a \in  \lie{qtr}_n$, we will 
denote the result of this action by $u \cdot a \in \lie{qtr}_n$. 

\begin{defn} The {\em $q$-divergence} is the map ${\rm qdiv}: \lie{tder}_n^{l\N}  \longrightarrow \lie{qtr}_n$ given by the formula 
\[ {\rm qdiv}(u = (a_1, \ldots a_n)) = \sum_{k=1}^n {\rm qtr} (x_k (\partial_k a_k))  \;\; \subset \;  \lie{qtr}_n . \] 
\end{defn} 

\begin{rem} 
{\rm If $l=2$ and $q=-1$, we recover the even tangential derivations, the supertrace  and the superdivergence.  Further on, 
we will consider the superdivergence as a particular case of $q$-divergence with $l=2$. } 
\end{rem} 

\begin{prop}  \label{qdiv} 
The $q$-divergence is a cocycle on $\lie{qtr}_n$, in the sense that 
$u \cdot {\rm qdiv}(v) - v \cdot {\rm qdiv}(u) = {\rm qdiv}([u,v])$ for all $u, v \in \lie{tder}_n^{l\N}$. 
\end{prop} 
{\it Proof:} Let $u = (a_1, \ldots, a_n)$, $v = (b_1, \ldots, b_n)$ such that the lengths of $a_i$ and $b_i$ are divisible by $l$ for all $i = 1, \ldots n$. Using the computation in the proof of Proposition 3.6 in \cite{AT}, we obtain: 
\begin{eqnarray*}
 {\rm qdiv} ([u,v])  & = &  \sum_{k=1}^n {\rm qtr} \{ x_k  u (\partial_k b_k) - x_k v (\partial_k a_k) \} +  \\ 
         &  + &   \sum_{k=1}^n {\rm qtr} \{  - x_k (\partial_k b_k) a_k  + x_k a_k (\partial_k b_k) +  x_k (\partial_k a_k) b_k  - x_k b_k (\partial_k a_k)  \} + \\ 
       & + &   \sum_{k,i=1}^n {\rm qtr} \{ x_k (\partial_i b_k)x_i (\partial_k a_i)   -  x_k (\partial_i a_k) x_i (\partial_k b_i) \}   .
       \end{eqnarray*} 
Consider the second line. Since the lengths of $a_k$ and $b_k$ are divisible by $l$, we have ${\rm qtr} ( - x_k (\partial_k b_k) a_k ) = 
 { \rm qtr} (- a_k x_k (\partial_k b_k)) $ and  ${\rm qtr} (  x_k (\partial_k a_k) b_k ) = 
 { \rm qtr} ( b_k x_k (\partial_k a_k)) $. The third line is zero, since the length of the expression $x_i (\partial_k a_i)$ 
 is divisible by $l$. Therefore, we have 
 \begin{eqnarray*} 
   {\rm qdiv} ([u,v]) & = &   \sum_{k=1}^n {\rm qtr} \{ x_k  u (\partial_k b_k) - x_k v (\partial_k a_k) + [x_k, a_k] (\partial_k b_k) -
 [x_k, b_k] (\partial_k a_k) \} =  \\ 
       & = &  \sum_{k=1}^n {\rm qtr}  \{ u(x_k(\partial_k b_k) ) - v(x_k (\partial_k a_k)) \}  = \\
       & = & u \cdot {\rm qdiv} (v) - v \cdot {\rm qdiv} (u) . 
\end{eqnarray*}  \qed

\section{The Grothendieck-Teichmuller Lie algebra $\lie{grt}_1$ and the Kashiwara-Vergne Lie algebra $\lie{krv}$.} 

\subsection{The Grothendieck-Teichmuller Lie algebra $\lie{grt}_1$.} 
Recall the Grothendieck-Teichmuller Lie algebra $\lie{grt}_1$ first introduced in \cite{Dr}. First we need to define the 
$n$-strand braid Lie algebra $\lie{t}_n$. 

\begin{defn} 
The Lie algebra $\lie{t}_n$ is generated by $n(n-1)/2$ elements $t^{i,j} = t^{j.i}$, where $1 \leq i,j \leq n$ and relations 
\[ [ t^{i.j}, t^{k,l} ] = 0 \;\;\;\; {\rm if} \;\;\;\; \{i,j\} \neq \{k,l \} ,  \;\;\; {\rm and}  \] 
\[ [ t^{i,j} + t^{i,k}, t^{j,k}  ] = 0, \]
for all triples of distinct indices $i,j ,k$. 
\end{defn} 

\begin{defn} 
 The Lie algebra $\lie{grt}_1$ is spanned by the elements $(0, \psi) \in \lie{tder}_2$, that satisfy the following  
relations: 
\begin{equation} \psi(x,y) = -\psi(y,x) ,  \label{asym} \end{equation} 
\begin{equation} \psi(x,y) + \psi(y,z) + \psi(z,x) =0 \;\;\;\; {\rm for} \;\;\;\; x+y+z=0,   \label{tri} \end{equation} 
\begin{equation} \psi(t^{1,2} , t^{2, 34}) + \psi(t^{12, 3}, t^{3,4}) = 
\psi(t^{2,3}, t^{3,4}) + \psi(t^{1,23}, t^{23,4}) + 
\psi(t^{1,2}, t^{2,3}),  \label{pent} 
\end{equation} 
where the last equation takes values in the Lie algebra $\lie{t}$, and  $t^{1,23} = t^{1,2} + t^{1,3}$, etc.   
\end{defn} 

\begin{exm} Let  $ \psi = [x, [x,y]] - [y,[y,x]] $. Then $(0, \psi) \in \lie{grt}_1$. 
\end{exm}  

The Lie algebra structure on $\lie{grt}_1$ is given by the Ihara bracket:
\[ \{\psi_1, \psi_2 \} = (0, \psi_1) (\psi_2) - (0, \psi_2) (\psi_1) + [\psi_1, \psi_2 ] . \] 

In \cite{AT}, the following statement was shown (Theorem 4.1): 

\begin{thm}The map $ \nu : \psi \longrightarrow (\psi(-x-y, x), \psi(-x-y, y))$ is an injective Lie algebra homomorphism  
 from $\lie{grt}_1$ to $\lie{tder}_2$. 
 \end{thm}

\subsection{Soul\'e elements} 

Consider the free Lie algebra 
\[ \lie{lie}_2 = \prod_{m=1}^\infty \lie{lie}^m (x,y) . \] 
Let ${\rm deg}( x) = {\rm deg}( y) = 1$. Then $\lie{lie}^m(x,y)$ is the graded component of weight $m$. This graded component  can be decomposed further into a direct sum 
of $\K$-vector spaces 
\[   \lie{lie}^m(x,y) = \oplus_{i=1}^{m-1} \lie{lie}^{(i, m-i)} (x,y) , \] 
where the subspace $ \lie{lie}^{(i, m-i)} (x,y) $ is spanned by the elements with $x$-degree $i$ and $y$-degree $(m-i)$. 

Similar decompositions exist for the universal enveloping algebra 
${\rm Ass}_2$ with   ${\rm deg}(x) = {\rm deg}(y)=1$. Then we have  
\[ {\rm Ass}_2 = \prod_{m=0}^\infty {\rm Ass}^m  (x, y) , \;\;\;\;\;\;\;\;\;\;\;   {\rm Ass}^m  (x, y)  = \oplus_{i=0}^m {\rm Ass}^{(i, m-i)} (x,y) , \] 
where ${\rm Ass}^{(i, m-i)}$ is spanned by the monomials with $x$-degree $i$ and $y$-degree $(m-i)$.  The bigrading factors 
through to the vector spaces ${\lie{ tr}}_2$, ${\lie {str}}_2$ and ${\lie {qtr}}_2$.  

\noindent {\bf Notation:} We will denote by $f^{(p,t)}$ the $(p,t)$-component of the element $f$ lying in any of the spaces 
$\lie{lie}_2$,  ${\rm Ass}_2$, $\lie{tr}_2$, $\lie{str}_2$ and $\lie{qtr}_2$. For example, if $f \in \lie{tr}_2$, then 
$f^{(p,t)} = f \, \cap \, \lie{tr}_2^{(p,t)}$.  


The following facts are known about $\lie{grt}_1$ (see, for example, \cite{Ih1}) : 
\begin{thm}   \label{soule} 
\begin{enumerate} 
\item  The Lie algebra $\lie{grt}_1 \ni (0, \psi)$ is graded by {\em weight} of the elements $\psi \in \lie{lie}_2$ with ${\rm deg} (x) = {\rm deg (y) } = 1$. The  first 
nontrivial component of $\lie{grt}_1$ lies in weight $m=3$. 
\[ \lie{grt}_1 = \oplus_{m=3}^\infty \lie{grt}_1^{m} .\]
It also admits a decreasing filtration by {\em depth}, the smallest $y$-degree in any monomial in $\psi$. 
We will denote the associated graded component of depth $n$ in $\lie{grt}_1$ by $\lie{grt}_1^{(n)}$. 
\item  For each odd integer $m \geq 3$ there exists a nonzero element $\sigma_m \in \lie{lie}^m (x,y)$, called the Soul\'e element,  
such that $(0, \sigma_m) \in \lie{grt}_1^m$. 
Hence the dimension ${\rm dim} (\lie{grt}_1^{m} ) \geq 1$ for any odd $m \geq 3$. 
\item   For each odd integer $m \geq 3$, the element $\sigma_m$ has depth $1$. The $\lie{lie}^{(m-1,1)}(x,y)$-component of $\sigma_m$ is proportional to $\ad_x^{m-1} y$  with a nonzero coefficient. 
\end{enumerate} 
\end{thm} 

\begin{exm} \cite{Ih2}  {\rm The Soul\'e elements below are normalized to have integer coefficients. }  \label{si3si5} 
\[ \begin{array}{ll} 
 \sigma_3 =  & [x, [x, y]] - [y,[y,x]] .\\
\sigma_5 =  &  2[x,[x,[x,[x,y]]]]  - 2 [y,[y,[y,[y,x]]]] +4[x,[x,[y,[x,y]]]] - 4[y,[y,[x,[y,x]]]] - \\ 
 &  -3[[x,y],[x,[x,y]]] + 3[[y,x],[y,[y,x]]] .  \end{array} \] 
\end{exm} 

\subsection{Kashiwara-Vergne Lie algebra $\lie{krv}$.} 

The Lie algebra $\lie{krv}$ was introduced in \cite{AT} as the Lie algebra of the group that acts freely and 
transitively on the space of solutions of the {\it Kashiwara-Vergne problem} \cite{KV}. 
(In \cite{AT}, the algebra $\lie{ krv}$ was denoted $\widehat{\lie{kv}_2}$).  

\begin{defn} 
\[ \lie{krv} := \{  u = (a,b) \in \tder_2  \; :\; \] 
\begin{equation}  [x,a] + [y,b] = 0,     \label{krv_eq1}   \end{equation} 
\begin{equation} \ddiv(u) = \tr( -f(x+y) +f(x)+f(y) ) \}   \label{krv_eq2}  \end{equation} 
for some element $f \in \tr_1$. 
\end{defn} 

\begin{exm} {\rm The element $t = (y,x)$ belongs to $\lie{krv}$ with $f =0$. } 
\end{exm} 

Then $\krv$ is a Lie subalgebra of the algebra of tangential derivations $\tder_2$ (see \cite{AT}, p. 14). 
Moreover, this Lie algebra contains the image of $\grt_1$ under the map $\nu$. 

\begin{thm}  The map $ \nu : \psi \longrightarrow (\psi(-x-y, x), \psi(-x-y, y))$ is an injective Lie algebra homomorphism  
 from $\lie{grt}_1$ to $\krv$. 
 \end{thm} 
The Lie algebra $\krv$ is known to contain a one-dimensional central subalgebra generated by the 
element $t= (y,x)$, and a Lie subalgebra isomorphic to $\grt_1$. It was conjectured in \cite{AT} that 
in fact $\krv  \cong \K t \oplus \grt_1$. One of the results of this paper, given below in Theorem \ref{divthm} confirms
 this conjecture up to the elements of depth $\geq 3$ in $\grt_1$.  
 
 To prove this theorem, we will need another  presentation of $\krv$. 
 
\begin{prop} There is a Lie algebra isomorphism between $\lie{krv}$  \label{krvprime}
 and the Lie subalgebra $\lie{krv'}$ of $\lie{tder}_2$ defined as follows: 
\[ \{ (0, \psi) \in \tder_2(z,y)  \; :  \,\, [y, \psi(z, y)] \in {\rm im}({\rm ad}_{y+z}), \] 
\[  {\rm tr}( y \partial_y \psi_{\geq 2}) = \tr( -f(-z) +f(-y-z) + f(y) ) \} , \] 
where $\psi=\psi_1+\psi_{\geq 2}$ with $\psi_1= c(z+y)$ of degree one and $\psi_{\geq 2}$ of degree at least two.
 \end{prop} 
{\it Proof :} The change of 
 variables $y=y, z=-x-y$ in $\mathfrak{lie}_2$ induces a map taking the derivation
\[ x \mapsto [x, a(x,y)], \hskip 0.5cm  y \mapsto [y, b(x,y)]  \] 
to the derivation
\[ z \mapsto 0, \hskip 0.5cm y \mapsto[y, b(-y-z, y)]. \] 
Let $\psi(z,y) = b(-y -z, y)$. 
Then the first condition (\ref{krv_eq1}) defining the algebra $\lie{krv}$ is satisfied if and only if $[y, \psi(z, y)] \in {\rm im}({\rm ad}_{y+z})$.\\
Now let $u=(a(x,y),b(x,y))$ be a solution to the equation (\ref{krv_eq1}). It implies \[\partial_x( [x,a(x,y)] + [y,b(x,y)])=0,\] which yields  \[x\partial_x a(x,y)= a(x,y) - y\partial_xb(x,y).\]
To check condition (\ref{krv_eq2}), we have to compute the divergence: 
\[ {\rm div}(u)(x,y)= {\rm tr}(x\partial_x a(x,y) + y \partial_y b(x,y)) = {\rm tr} (a(x,y)) + {\rm tr} (y (\partial_y b(x,y) - \partial_x b(x,y))) \]
Changing variables we obtain 
\begin{align*} {\rm div}(u)(-z-y,y)=& {\rm tr} (a(-z-y,y)) + {\rm tr} (y (\partial_y b(-z-y,y) - \partial_x b(-z-y,y))) \\
						=&{\rm tr} (a(-z-y,y)) + {\rm tr} (y \partial_y \psi(z,y)).\end{align*}
The trace ${\rm tr} (a(-z-y,y))$ is zero on all elements $a \in \lie{lie}(z,y)$ of degree at least two, but not on elements of degree one. Writing $\psi$ as a sum $\psi = \psi_1+\psi_{\geq 2}$, one can easily show that $\psi_1$ is proportional to $(z+y)$, and 
$\psi_{\geq 2}$ has to satisfy the condition 
$\tr ( y \partial_y \psi_{\geq 2}) = \tr( -f(-z) +f(-y-z) + f(y) )$. 

\qed 
 

Morally, the algebra $\lie{krv'}$ is the same algebra $\lie{krv}$ viewed in a different coordinate system, obtained by 
a nonsingular linear change of variables. 

\begin{cor}  The map $\psi(x,y) \longrightarrow \psi(z,y)$   is a an injective homomorphism of Lie algebras 
$\grt_1  \longrightarrow \krv'$      \label{grt_1_to_krv'}. 
\end{cor} 

{\it Proof:}   The map is the composition of the homomorphism 
\[ \nu : \psi (x,y) \longrightarrow (\psi(-x-y, x), \psi(-x-y, y)) = 
(a(x,y), b(x,y)) \in \lie{krv} \] 
and the isomorphism $\lie{krv} \cong \lie{krv'}$ of Proposition \ref{krvprime} that maps 
\[ (a(x,y), b(x,y)) \longrightarrow (0, b(-y-z, y)) = (0, \psi(z,y) ). \]    
\qed

\subsection{The target space of the injection of $\lie{grt}_1$ into $\lie{krv'}$ in depth $2$.}

Similarly to the algebra $\lie{grt}_1$, $\lie{krv'}$ admits a grading by weight with $\on{deg}(y) = \on{deg}(z) =1$ and 
a filtration by depth. 
The depth filtration in $\lie{krv'}$ 
is defined as follows: $\psi$ is of depth $n$ (denoted $\psi \in \krv'^{(n)}$) if all Lie words in $\psi(z,y) $ contain at least $n$ $y$'s. 
The $n$-th graded component of the depth-associated graded Lie algebra is defined as
\[ \textbf{gr} (\krv')^{(n)} = \krv'^{(n)}/\krv'^{(n+1)}. \] 

It is clear from the construction that the injective homomorphism $\grt_1 \longrightarrow \krv'$  of Corollary \ref{grt_1_to_krv'}
is compatible with the filtration by depth. \\
To see if the injective homomorphism  $\grt_1 \longrightarrow \krv'$ is also surjective in depth $2$, our strategy will be to compute 
an upper bound for the dimension of the graded component  $\bf{gr}(\lie{krv'})^{(2)}$, and compare it with the known lower bound for the dimension of 
${\bf gr}(\lie{grt}_1)^{(2)}$.

To get upper bounds of the dimensions of the associated graded components $ \textbf{gr} (\krv')^{(n)}$,  it is enough to consider the graded vector space $\overline{\krv}$ that is 
obtained by dropping the terms of depth greater than $1$ in the relations for the $\krv'$.  Then we have 

\[ \overline{\krv} = \{ (0, \psi) \in \tder_2(z,y)  \; :  \; {\rm tr}( y \partial_y \psi) = \tr(y g(z)) \}, \]
where $g(z)$ is  a formal power series in $z$. Indeed, for the elements $\psi$ of depth greater than $1$, the relation  
${\rm tr}( y \partial_y \psi) = \tr( -f(-z) +f(-y-z) + f(y) ) $ implies ${\rm tr}( y \partial_y \psi) =0$, and for the elements of depth $1$ we obtain ${\rm tr}( y \partial_y \psi) =  \tr(y g(z))$.  The relation $[y, \psi(z, y)] \in {\rm im}({\rm ad}_{y+z})$ is trivial 
for $\psi$ of any given depth $n$, because in this case it reads 
\[ [y, \psi^{(n)}(z,y)] = {\rm ad}_y \phi^{(n)}(z,y) + {\rm ad}_z \phi^{(n+1)}(z,y),  \]  
where $\phi^{(n)}$ and $\phi^{(n+1)}$ are some elements of depth $n$ and $n+1$  respectively. Such elements always 
exist (for instance, one can take $\phi^{(n)} = \psi^{(n)}$ and $\phi^{(n+1)} =0$.) 

The  graded vector space $\overline{\krv}$ 
 admits 
 a filtration by depth, the smallest $y$-degree in any monomial of a given expression. 
There is a natural injection of $\textbf{gr}(\krv') \to \overline{\krv}$: for $\psi$ of depth $n$ we simply take its part containing exactly $n$ 
$y$'s. Hence, the dimensions of the graded components of $\overline{\krv}$ in each depth give upper bounds for the dimensions of the graded components of $\krv'$ in the same depth.  In particular, the dimensions of the weight-graded components of 
$\overline{\krv}$ of depth $2$ give upper bounds for the dimensions of the corresponding 
components of $\textbf{gr}(\krv')^{(2)}$.  
Note that for $\psi \in \overline{\krv}$ of depth $2$, the relation 
${\rm tr}( y \partial_y \psi) = \tr(y g(x)) $ 
implies $g(x) =0$.  Therefore, to estimate the dimensions of the graded components  of ${\textbf gr}(\krv')^{(2)}$, 
we need to find the dimensions of the graded components of the space 
\[ \textbf{gr}(\overline{\krv})^{(2)} = \{ \psi \in \lie{lie}(x,y)^{(2)} /\lie{lie}(x,y)^{(3)} \; : \; {\rm tr} (y \partial_y \psi) = 0 \} . \] 
This computation  will be carried out in Section 5, as a particular case of a more general computation.

\section{The action of the  $q$-divergence cocycle on $\nu(\lie{grt}_1)$.} 

Recall that for $q^l =1$, with $l \geq 2$, the $q$-divergence cocycle is defined on the elements of $\lie{tder}_n$ 
with length of words divisible by $l$, by the formula $\qdiv(u = (a_1, \ldots a_n)) = \sum_{k=1}^n \qtr (x_k (\partial_k a_k))$.  
In particular, the case $l=1$ corresponds to divergence, and $l=2$ to the superdivergence. 
We would like to consider the action of this cocycle on $\nu(\lie{grt}_1) \subset \lie{\tder}_2$. 

\subsection{The $q$-divergence cocycle on the elements of depth $2$.}

 \begin{exm} {\rm  Let $\psi = \sigma_3 = [x, [x,y]] - [y,[y,x]]  \in \lie{grt}_1$. 
 It is easy to find that  $\nu (\psi) = (  [y, [y,x]],    [x, [x, y]])$. 
Then we have, ${\rm div}(\nu (\psi) ) = {\rm tr} ( x y^2 + y x^2) \neq 0$. The superdivergence cannot be applied since the 
length of $\nu(\psi)$ is odd. If $q^3=1$, then the $q$-divergence can be calculated: 
$ {\rm qdiv}(\psi) = {\rm qtr}(x y^2 + y x^2) \neq 0$. } 
\end{exm} 

However, we know that the divergence cocycle has a large kernel on $\nu(\lie{grt}_1)$: 
 
\begin{prop}  (\cite{AT}, page 14 and Theorem 4.1).  \\
\noindent   The commutant $\{ \lie{grt}_1, \lie{grt}_1 \}$ lies in the kernel of the map 
\[ {\rm div} \circ \nu : \lie{grt}_1 \longrightarrow \lie{tr}_2 . \] 
\end{prop} 

We would like to find if, similarly, the $\qdiv$ cocycle for $q$ a root of unity of degree $l \geq 2$  can have 
nontrivial kernel on $\nu(\grt_1)$. It is not hard 
to see that, in contrast with $\ddiv$, the cocycle $\qdiv$ for $l \geq 2$ is nonzero on all commutators of the Soul\'e 
elements: 
\[   \qdiv (\nu ( \{ \sigma_m , \sigma_n \} )) \neq 0 . \] 
Here $m<n$ are odd numbers, $m,n \geq 3$ and in case of $q$-divergence, $m+n$ is divisible by $l$. 
For example, Lemma 2.5 in \cite{Ma} contains an explicit formula for $\{ \sigma_m , \sigma_n \}$ modulo  the 
elements of depth $3$, that allows to prove the claim by a direct computation.  
Below we will aim to  obtain a more general result: if $q \neq 1$, the cocycle $\qdiv$ has a trivial kernel 
on ${\bf gr}(\nu(\grt_1))^{(2)}$.

First we will derive a simple formula for the action of  $\qdiv$, whenever applicable, on the elements of depth $2$ in 
$\nu(\lie{grt}_1)$. This formula is also applicable in case of $\ddiv$.

\begin{lem}  \label{nupsi}
Suppose $(0,\psi) \in \lie{grt}_1$ is of weight $2n$ and of depth $2k$ for some natural $k <n$. 
Then for $q$ a primitive $l$th root of unity with $l \geq 1$, whenever the $q$-derivation can be applied, we have 
\[ {\rm qdiv} (\nu (\psi ))^{(2n-2k, 2k)}  = {\rm qtr} \left( y \partial_y  \left(\psi(x,y)^{(2n-2k, 2k)} \right)\right) . \] 
\end{lem} 

{\it Proof:}  Assume that 
the weight $2n$ of $\psi$ is divisible by the order of $q$, so that the $q$-derivation can be applied. Then we have 
\[ {\qdiv} (\nu(\psi (x,y) ) ) = {\qtr} \left( x \partial_x \psi (-x-y,  x) + y  \partial_y \psi(-x-y, y) \right) . \]  
Using the equations (\ref{tri}) and (\ref{asym}), we can write 
\[ \psi (-x-y, y) = \psi(-x-y, x) + \psi(x,y) . \] 
Then 
\[ {\rm qdiv} (\nu(\psi (x,y) ) ) = {\rm qtr} \left( x \partial_x \psi (-x-y,  x) + y  \partial_y \psi(-x-y, x)  + y \partial_y \psi(x,y) \right) . \]  
The expression $x \partial_x \psi (-x-y,  x) + y  \partial_y \psi(-x-y, x)$ differs from $\psi(-x-y, x)$ by removing the last letter of  
each monomial and inserting it in front. Under the respective trace maps, it gives 
 \[ \begin{array}{ll} 
 {\rm qtr} \left( x \partial_x \psi (-x-y,  x) + y  \partial_y \psi(-x-y, x) \right)  = & q \cdot {\rm qtr}( \psi(-x-y, x)). 
\end{array} \]

Now we will use the assumptions on the weight and depth of $\psi$ to analyze the expression  ${\qtr}(\psi(-x-y, x))$.   Using the equation (\ref{tri}) again, we can write 
\[  {\qtr} (\psi(-x-y, x)) = {\qtr} ( \psi(-x-y, y) - \psi(x, y) ). \] 

We only need to consider the lowest possible $y$-degree component of $( \psi(-x-y, y) - \psi(x, y) )$.
 Since $\psi(x,y)$ has no component 
with $y$-degree less than $2k$,  the minimal possible $y$-degree in $\psi(-x-y, y)$ is also $2k$, and the component of 
$\psi(-x-y, y)$ with the minimal $y$-degree is obtained by choosing the summand 
with $-x$ each time $(-x-y)$ occurs in the expression. We have   
\[ \psi(-x-y, y)^{(2n-2k, 2k)} = \psi(-x, y)^{(2n-2k, 2k)} = 
(-1)^{2n-2k} \psi(x, y)^{(2n-2k, 2k)} = \psi(x, y)^{(2n-2k, 2k)}. \] 
Then we obtain 
\[  {\qtr} \left(\psi(-x-y, x)^{(2n-2k, 2k )}  \right) = {\qtr} \left( (\psi(x, y) - \psi(x, y))^{(2n-2k, 2k)} \right) = 0. \] 
Finally, 
\[ {\qdiv}  (\nu(\psi (x,y)  ) )^{(2n-2k, 2k)}  = {\qtr} \left( y  \partial_y ( \psi(x, y)^{(2n-2k, 2k)} \right) . \] 
\qed  

 It is well known that the weight-associated 
graded component ${\bf gr}( \grt_1)^{(2)}$ is spanned by the corresponding parts of the 
commutators of the Soul\'e elements. All these elements have even total weight, and $y$-degree equal to $2$. Then  
by Lemma \ref{nupsi}, studying the action of $\qdiv$ on this space is equivalent to studying the action of $\qtr (y \partial_y \cdot)$  on the elements of $\lie{lie}(x,y)$ of even weight and $y$-degree $2$. 
 

Therefore, we are interested in the dimensions of the weight graded components of the space 
\[ \on{ker}(\qdiv)({\bf gr}(\nu(\lie{grt}_1))^{(2)}) = \{ \psi \in \lie{lie}(x,y)^{(2)} /\lie{lie}(x,y)^{(3)} \; : \; {\qtr} (y \partial_y \psi) = 0 \} . \] 
This computation will be given in the next section.

\section{The upper bound for $\bf{gr}(\lie{krv'})^{(2)}$  and the kernel of $q$-divergence on ${\bf gr}(\nu(\lie{grt}_1))^{(2)}$.} 

The set of relations that determines the upper bound for the graded components of 
 $\bf{gr}(\lie{krv'})^{(2)}$ (Section 3), is a particular case (when $q=1$) of the relations that define the kernel of the $q$-divergence on the image of $\lie{grt}_1$ in depth $2$ (Section 4). We will use the same combinatorial method to establish the dimensions of the graded components in both cases.  

\subsection{Combinatorial reformulation of the question.} 

Let $q$ be a primitive $l$th root of unity with $l \geq 1$. We are interested in the dimensions of the graded components of 
 the space  
 \[  \{ \psi \in \lie{lie}(x,y)^{(2)} /\lie{lie}(x,y)^{(3)} \; : \; {\qtr} (y \partial_y \psi) = 0 \} . \] 

We are interested in the action of ${\qtr}( y \partial_y  \; \cdot) $ on the elements of the free Lie algebra 
in two variables $(x,y)$  with exactly $2$ $y$'s in any of its Lie words. In other words, we are interested in 
the kernel of the operator $ {\qtr}( y \partial_y  \; \cdot) $ acting on  the 
elements of the form 
\[ \sum_{k} c_{l,k} [ {\rm ad}_x^k y, {\rm ad}_x^l y ] , \]
where $c_{l,k} \in \K$ are some coefficients. Because the operator $ {\qtr}( y \partial_y  \; \cdot)$ 
preserves the weight, it is enough to consider its action on the graded components of a given weight $n$, 
\[ \sum_{k=0}^n c_k  [ {\rm ad}_x^k y, {\rm ad}_x^{n-k} y ] , \]
with some coefficients $c_k$, $ k = 0, \ldots, n$, for every natural $n$ divisible by $l$.

\begin{prop}   For each $n$ divisible by $l$, the dimension of the vector space  \label{qkern} 
\[ \on{qker}_n =   \{ \psi \in \lie{lie}(x,y), \;  \psi = \sum_{k=0}^n c_k  [ {\rm ad}_x^k y, {\rm ad}_x^{n-k} y ] , \;  \; {\qtr}( y \partial_y  \; \psi) =0 \} \]
is equal to the dimension of the vector space of homogeneous polynomials of degree $n$ in two variables $p(v, w)$ satisfying  
 \begin{align*}  p(v, w)  &= - q \cdot p(v + w, (q-1)v - w) \\
  p(v, w) &= - p(w, v) 
  \end{align*}

In particular, for $q=1$, we obtain the following relations: 
 \begin{align*}  p(v, w)  &= - p(v + w,  - w) \\
  p(v, w) &= - p(w, v) 
  \end{align*}
\end{prop}

{\it Proof:} 
The generating series for the expressions of the form $\sum_{k=0}^n c_k  [ {\rm ad}_x^k y, {\rm ad}_x^{n-k} y ]$  for all natural $n$ can be written as    
\[ \phi(\alpha, \beta) = [ e^{\alpha \; {\rm ad}_x }y , e^{\beta \; {\rm ad}_x } y] . \] 
 Consider the action of ${\qtr}( y \partial_y  \; \cdot) $ on this generating series :
\begin{flalign*}
{\qtr}( y \partial_y \; \phi(\alpha, \beta) ) & = 
 {\qtr}( y \partial y (e^{\alpha x} y e^{- \alpha x} e^{\beta x} y e^{-\beta x} - e^{\beta x} y e^{- \beta x} e^{\alpha x} y e^{-\alpha x})) \\
&={\qtr}( y (e^{\alpha x} y e^{- \alpha x} e^{\beta x} - e^{\beta x} y e^{- \beta x} e^{\alpha x} ))\\
&= q \cdot  {\qtr} (e^{\alpha x} y e^{- \alpha x} e^{\beta x} y )
- q \cdot {\qtr}(e^{\beta x} y e^{- \beta x} e^{\alpha x} y ) .
\end{flalign*}
Here we have used $e^{-\beta x}=(\sum_{n=1}^\infty \frac{(-\beta)^m x^m} {m!})$ which implies that $e^{\alpha x} y e^{- \alpha x} e^{\beta x} y e^{-\beta x}$ only has one term ending with $y$. 
Thus ${\qtr}( y \partial_y \; \phi(\alpha, \beta) )$  is  an infinite series with terms of the form ${\qtr}( x^l y x^m y )$, 
that we have to consider up to the $\qtr$-symmetry: 
\[ \qtr(x^m y x^k y) = q^{(m+1)} \qtr (x^k y x^m y) = q \cdot \qtr ( x^k y (q x)^m y) .  \]
These expressions are in one-to-one correspondence with the sums of monomials in two commuting variables, 
\[ u^m v^k  + q \cdot (qv)^m u^k  . \] 
Applying this correspondence, we transform $\qtr (y \partial_y \phi(x,y))$  into 
\[ T(\alpha, \beta; u,v) \equiv q \cdot  ( e^{\alpha u} e^{ (\beta - \alpha) v} ) -q \cdot (e^{\beta u} e^{(\alpha - \beta)v} ) + 
q^2 \cdot ( e^{\alpha  q v} e^{ (\beta - \alpha) u} )  - q^2 (e^{\beta q v} e^{(\alpha - \beta) u} ) .\] 
 Differentiating $T(\alpha, \beta; u,v)$ several times with respect to $\alpha$ and $\beta$ and setting these  parameters 
 equal to zero allows to recover the part of the expression of a given degree in $u$ and $v$. In particular, 
 if we want to recover the action of ${\qtr}( y \partial_y \cdot)$ on 
 $\psi = \sum_{k=0}^n c_k  [ {\rm ad}_x^k y, {\rm ad}_x^{n-k} y]$, 
 we have to apply the differential operator $p\left( \frac{\partial}{\partial \alpha}, \frac{\partial}{\partial \beta} \right) \vert_{\alpha=0, \beta=0}$,  where $p(t,s) = \sum_{k=0}^n c_k s^k t^{n-k}$. Applying $p\left( \frac{\partial}{\partial \alpha}, \frac{\partial}{\partial \beta} \right) \vert_{\alpha=0, \beta=0}$ to $T(\alpha, \beta; u,v)$ gives

\[ q \cdot p(u-v, v) -q \cdot p(v, u-v) + q^2 \cdot p(qv-u, u) - q^2 \cdot p(u, qv -u) . \] 
Because of the  Lie antisymmetry of $\psi$, 
the polynomial $p(t,s)$ satisfies the condition $p(t,s) = - p(s,t)$ with respect to the swap of the variables.
Therefore, 
we have 
\[ \left . p\left( \frac{\partial}{\partial \alpha}, \frac{\partial}{\partial \beta} \right) \right|_{\alpha=0, \beta=0} \, (T(\alpha, \beta; u,v) = -2 q \cdot  p(v, u-v) -2 q^2 \cdot  p(u,qv-u) .\] 
Setting this expression equal to zero, we have the equation 
\[ p(v, u-v) = - q \cdot p(u, qv-u) .\] 
With the change of variables $v \to v$, $u-v \to w$, we obtain the condition 
\[  p(v, w) = -q \cdot p(v+w, (q-1)v -w) , \] 
 which together with the antisymmetry yields the system announced in the Proposition. 
 \qed

\subsection{Upper bounds for the dimensions of the graded components of ${\bf{gr}}(\lie{krv'})^{(2)}$.} 

For each $n \geq 1$, denote by $\on{ker}_n$ the following vector space 
\[ \on{ker}_n =   \{ \psi \in \lie{lie}(x,y), \;  \psi = \sum_{k=0}^n c_k  [ {\rm ad}_x^k y, {\rm ad}_x^{n-k} y ] , \;  \; {\tr}( y \partial_y  \; \psi) =0 \} \]
 Then 
by Proposition \ref{qkern}, the dimension of the space $\on{ker}_n$ 
is equal to the dimension of the space of homogeneous polynomials of degree $n$ in two variables, satisfying the relations 
 \begin{align*}  p(v, w)  &= - p(v + w,  - w) \\
  p(v, w) &= - p(w, v) 
  \end{align*}

According to the discussion at the end of Section 3, the dimensions of $\on{ker}_n$ are equal to the dimensions of 
the weight associated graded components of  
${\bf{gr}}(\overline{\lie{krv}})^{(2)}= \{ \psi \in \lie{lie}(x,y)^{(2)} /\lie{lie}(x,y)^{(3)} \; : \; {\rm tr} (y \partial_y \psi) = 0 \}$
and 
give the upper bounds for the 
dimensions of the weight associated components of ${\bf{gr}}(\lie{krv'})^{(2)}$. We will compute these dimensions.

Note that the given transformations of homogeneous polynomials are implemented by 
the action of the group generated by the following  changes of variables:     

\[  a = \begin{pmatrix} 0 &1 \\ 1 &0 \end{pmatrix}, \quad   \begin{pmatrix} 0 &1 \\ 1 &0 \end{pmatrix}  \begin{pmatrix} v \\ w \end{pmatrix} = \begin{pmatrix}  w \\ v \end{pmatrix} \]
and
\[ b =  \begin{pmatrix*}[r] 1 &1 \\ 0 & -1 \end{pmatrix*}, \quad  \begin{pmatrix*}[r] 1 &1 \\ 0 & -1 \end{pmatrix*}   \begin{pmatrix} v \\ w \end{pmatrix} = \begin{pmatrix}  w + v \\ -w\end{pmatrix}  \]
This  is the dihedral group $D_{12} = \{s^2=1, r^6=1, srs^{-1}=r^{-1} \}$, 
where $s= a$ and $r=  b a$.  

This group acts on polynomials by $g \cdot p(v,w) = p(g^{-1}(v, w))$. Since the action is homogeneous, it 
 induces a representation of $D_{12}$ on the $(n+1)$-dimensional space of coefficients of degree $n$ homogeneous polynomials in two variables. 
Then the dimension of $\on{ker}_n$ for each natural $n$ is equal to the dimension of the subrepresentation where 
both $a=s$ and $b = rs$ act as $-1$. 
 
 \begin{exm} {\rm Let us consider the action of $s$ and $rs \in D_{12}$ on the homogeneous polynomials of degree $n=1$, 
$p(v,w) = a v +bw$. We have $s \cdot p(v,w) = p(s(v,w)) = p(w,v) = b v + a w$, and $rs \cdot p(v,w) = p((rs)^{-1} (v,w)) = 
p(v+w, -w) = av + (a-b)w$. For both $s$ and $rs$ to act by $-1$, the conditions on the coefficients $(a,b)$ are 
$b=-a, a=-b, a=-a, (a-b)=-b$ which yields only the trivial solution. Hence the dimension of $\on{ker}_1$ is zero. } 
\end{exm}  
   
  In general, our strategy will be to decompose the action of $D_{12}$ on the degree $n$ homogeneous polynomials 
  in two variables for each natural $n$ into a direct sum of irreducible 
representations, and find the dimension of the subspace where both elements $a =s$ and $b = rs$ act as $-1$. 

\begin{thm} The dimension of $\on{ker}_n$ is as follows:     \label{dim_kern} 
\[  \on{dim} \, (\on{ker}_n) = \left[ \begin{array}{ll} 
                                                0, & n \; {\rm odd} ; \\ 
                                                \left[ \frac{n}{6} \right]  &  n \; {\rm even}  
                                                \end{array} \right .   \]  
\end{thm} 
{\it Proof:} 
Here is the character table of $D_{12}$, where we denoted   by $\chi_{11}, \chi_{12}, \chi_{13}, \chi_{14}$ the caracters 
of the irreducible representations of degree $1$, and by 
 $\chi_{21}$ and $\chi_{22}$ those of degree $2$: 
 \[ \begin{array}{c|c|c|c|c|c|c}
   & e  & s & r & r^2 & r^3 & rs \\ \hline
   \chi_{11} & 1  & 1 & 1 & 1 &  1 &   1 \\
  \chi_{12}  & 1 &  1 & -1& 1 & -1 & -1 \\
   \chi_{13} & 1& -1 & 1 & 1 & 1 & -1 \\
   \chi_{14} & 1 & -1 & -1 & 1 & -1 & 1 \\
   \chi_{21} & 2 & 0 & 1 & -1& -2 & 0 \\
   \chi_{22} & 2 & 0 &-1 & -1 &2 & 0  
\end{array}  \]

The only representation that satisfies our requirements is $\chi_{13}$, that has the property $\chi_{13}(s) =\chi_{13}(rs)=-1$. 
We will compute the character $\chi_n$ of the representation of $D_{12}$ on the $(n+1)$-dinensional space of the coefficients 
of the homogeneous polynomials in two variables, as defined above. Then for each $n$ we will 
look for the multiplicity of $\chi_{13}$ in $\chi_n$.  

Computation given in the Appendix leads to the following character table for $\chi_n$: 

\[ \begin{array}{c|c|c|r|r|c|c|c}
& e  & s & r & r^2 & r^3 & rs &  \\ \hline
     \chi & n+1  & 1 & 1 & 1 &  n+1 &  1 & \text{ if } n \equiv 0 \text{ mod } 6  \\ \hline
   \chi & n+1  & 0 & 1 & -1 &  -(n+1) &  0 & \text{ if } n \equiv 1 \text{ mod } 6  \\ \hline
      \chi & n+1  & 1 & 0 & 0 &  n+1 &  1 & \text{ if } n \equiv 2 \text{ mod } 6  \\ \hline
   \chi & n+1  & 0 & -1 & 1 &  -(n+1) &  0 & \text{ if } n \equiv 3 \text{ mod } 6  \\ \hline
    \chi & n+1  & 1 & -1 & -1 &  n+1 &  1 & \text{ if } n \equiv 4 \text{ mod } 6  \\  \hline
   \chi & n+1  & 0 & 0 & 0 &  -(n+1) &  0 & \text{ if } n \equiv 5 \text{ mod } 6  \\ 
   \end{array} \]

For each $n$ we use the orthogonality of characters to compute the multiplicity of $\psi_{13}$ in the decomposition. We get 
\[
\langle \chi, \chi_{13} \rangle = \frac{1}{12} \sum_{g \in D_{12}} \chi(g) \cdot \chi_{13}(g^{-1}) =
 \left \{ \begin{array}{ll} 
0  & \text{if} \;  n \;\; \text{odd} \\ 
  \frac{n}{6} & \text{if} \;  n \equiv 0 \pmod{6}  \\
 \frac{n-2}{6} & \text{if} \;n \equiv 2 \pmod{6}\\
  \frac{n-4}{6} & \text{if} \;n \equiv 4 \pmod{6}\\
  \end{array} \right.\]

Therefore, the multiplicity of the character $\chi_{13}$ in the character $\chi_n$ for even $n$ equals $\left[\frac{n}{6} \right]$. 
This is the dimension of the space of the degree $n$ homogeneous polynomials in two variables that satisfy the 
conditions $p(w+v, -w) = -p(v,w)$ and $p(w,v) = - p(v,w)$. Then applying Proposition \ref{qkern}, we deduce that 
 the dimension of $\on{ker}_n$ is $\left[\frac{n}{6} \right]$ for even $n$, and zero for odd $n$.  
\qed 

\begin{rem}
{\rm In  \cite{Z}, Zagier computed graded dimensions of the double shuffle Lie algebra in depth 2. According to 
\cite{IhO}, his calculation gives rise to a pair of functional equations
$$
p(u, w)=-p(w,u) \hskip 0.3cm , \hskip 0.3cm p(u+w,u)+p(u+w,w)=0.
$$
These equations do not coincide with ours, but they also lead to a representation of the group $D_{12}$, 
 isomorphic to the one we obtained above. }
\end{rem}

\subsection{Surjection $\lie{grt}_1 \longrightarrow \lie{krv'}$ in depth $2$.} 

\begin{thm} The injective Lie algebra homomorphism $\lie{grt}_1 \longrightarrow \lie{krv'}$ is surjective on   \label{divthm}
the elements of depth $2$ modulo the elements of higher depth. 
\end{thm} 

{\it Proof:} 
Theorem \ref{dim_kern} provides an upper bound on the dimensions of the associated graded components 
of depth $2$ in each degree of the algebra $\lie{krv'}$. The weight of an element of $\lie{krv'}$ (total number of $x$'s and 
$y$'s in a Lie word) is related to the degree $n$ as follows: $k=n+2$. Then for the elements of weight $k$, this upper bound is equal to  $\left[ \frac{k-2}{6} \right]$ for even $k$ and $0$ for odd $k$. 

On the other hand, we can compute a lower bound for the dimension of the associated graded 
components of depth $2$ in $\lie{grt}_1$ for each given weight. 
We know that all commutators of the Soul\'e elements $\{\sigma_i, \sigma_{k-i}\}$, $3 \leq i \leq k-3, i \; \text{odd}$,  taken modulo the elements of depth $3$, are among the depth $2$ elements of weight $k$. 
Taking into account the antisymmetry of the Ihara bracket, there are $\left[ \frac{k-4}{4} \right]$ different  commutators 
of Soul\'e elements in the depth-associated graded component of a weight $k$ and depth $2$. 
In addition, it is shown by Goncharov \cite{G}, and independently by Ihara and Takao in \cite{Ih2}, 
 that the dimension of the space of all linear relations 
of the form 
 \[\sum_{i=3 \; i \text{ odd}}^{k-3} a_i \{\sigma_i, \sigma_{k-i}\} \equiv 0 \; \pmod {\lie{grt}_1^{(3)}} \]  is given by :
\[ \left[ \frac{k-4}{4}\right] -\left[ \frac{k-2}{6} \right] . \]
Therefore, the dimension spanned by the commutators of the Soul\'e elements of  weight $k$ and depth $2$, modulo elements of higher depth, 
is exactly $\left[ \frac{k-2}{6}\right]$. This is the lower bound for the dimension of the depth-associated graded component 
of $\lie{grt}_1$ of depth $2$ and weight $k$. Since $\lie{grt}_1  \mapsto \lie{krv'}$ is a depth-preserving injection, and 
the lower bound for $\lie{grt}_1$ coincides with the upper bound for $\lie{krv'}$ in the depth-associated graded component of depth $2$ in each weight, 
we have a surjection on the elements of depth $2$, modulo the elements of higher depth. 
\qed

\begin{exm} {\rm The first degree where the dimension of $\on{ker}_n$ is positive is $n=6$, which corresponds 
to the  Lie words of degree $8$. Since we know that $\lie{grt}_1$ coincides with $\lie{krv'}$ in depth $2$, the 
depth $2$ graded commponent of $\lie{krv'}$ of degree $8$ is $1$-dimensional and is spanned by the element } 
  \[ \{\sigma_3, \sigma_5\} = 5 [ {\rm ad}_x^4 y, {\rm ad}_x^2 y ] + 2 [ {\rm ad}_x^5 y, {\rm ad}_x y ] .\]
\end{exm}  

\begin{rem}
{\rm L. Schneps explained to us the her results in \cite{S} imply that the equation ${\rm tr}(y \partial_y \psi)=0$ in depth 2 is equivalent to the stuffle equation in the definition of the double shuffle Lie algebra. This argument shows that the double shuffle Lie algebra and $\mathfrak{krv}'$ are isomorphic in depth 2. Hence, the calculation by Zagier \cite{Z} of graded dimensions of the double shuffle Lie algebra can be used to give an alternative proof of our result above.}

\end{rem}


\subsection{The kernel of the $q$-divergence cocycle on  $\nu(\lie{grt}_1)$ in depth $2$.} 

\begin{thm} Let $l$ be a positive integer, $l \geq 2$, and $q$ a primitive $l$th root of unity. 
The kernel of the $q$-divergence cocycle on the elements 
of ${\bf gr}(\nu(\lie{grt}_1))^{(2)}$ with weight divisible by $l$, is trivial.     \label{qker} 
\end{thm} 

In view of the discussion in Section 4, and Proposition \ref{qkern},  Theorem \ref{qker} will follow from  Proposition 
\ref{qker_poly} below.

\begin{prop} Let $q$ be a primitive $l$th root of unity with $l \geq 2$. 
No homogeneous polynomial in two variables satisfies the conditions  \label{qker_poly} 
\begin{align*}  p(v, w) &= -q \cdot p(v+w, (q-1)v -w) ,\\
  p(v, w) &= - p(w, v) .
  \end{align*}
\end{prop} 

{\it Proof:} 
To find the homogeneous polynomials in each  degree divisible by $l$ that satisfy the obtained conditions, consider the transformation 
of the variables  
\[  b :  \left( \begin{array}{c} v \\ w \end{array}  \right) \to \left( \begin{array}{c} v+w \\  (q-1)v-w \end{array} \right) ; \quad \quad 
a : \left(  \begin{array}{c} v \\ w \end{array} \right) \to \left( \begin{array}{c} w \\ v \end{array}  \right) . \] 

Then applying   $(ab)^{-1}$ to the polynomial $p(v,w)$, we obtain the following equation for an antisymmetric homogeneous 
polynomial: 
\[ p(v,w) = q \cdot (ab)^{-1} \cdot p(v,w) =  q \cdot p((ab) \cdot (v,w)) = q \cdot p((q-1)v-w , v+w) .\] 

The transformation of variables $ab$ is represented by the following matrix: 
\[ ab = \left( \begin{array}{cc} 
                    q-1 & -1 \\ 
                    1 & 1  \end{array} \right) . \] 
This matrix is                     
diagonalizable with  eigenvalues $\lambda_{1,2} = \frac{q \pm \sqrt{q^2-4q}}{2}$.  
The product of these eigenvalues is $\lambda_1 \cdot \lambda_2 = q$. Therefore, with a suitable linear 
change of variables in the polynomial $p(v,w)$, we have an equation for  another homogeneous polynomial $P(V,W)$ 
in variables $V,W$ of the same total degree as $p(v,w)$: 
\[ P(V,W) = q \cdot P(\lambda_1 V, \lambda_2 W) . \] 
To solve this equation, we have to satisfy the relations of the form 
\[  V^t W^s = q \cdot (\lambda_1)^t (\lambda_2)^s V^t W^s \] 
for each monomial in $P(V,W)$. Suppose that $t \geq s$ (the case $t \leq s$ is solved by the same argument), 
then we have  
\[ 1 = q \cdot (\lambda_1 \lambda_2)^{s} \lambda_1^{t-s}  =  q \cdot q^s \lambda_1^{t-s} . \] 
If $t  >s$, then  the absolute value of $\lambda_1$ (and therefore also of $\lambda_2$) has to be $1$.  
Consider  the equation 
\[ \lambda = \frac{q \pm \sqrt{q^2-4q}}{2}  = e^{i \phi} \] 
where $\phi$ is a real number. The only solution of this equation is  $\phi = \pm \frac{\pi}{3}$ and $q =1$.  
This corresponds to the case of divergence, that we considered before, and where obtained many nontrivial solutions.  
 
For the $q$-divergence with $l \geq 2$, we have to have $t=s$, and moreover $1 = q^{s+1}$. 
By the assumption on the degree of the polynomial $p(v,w)$, we have $s+t +2 = 2s+2 = lk$ for some positive integer $k$. Then $q^{s+1}= q^{\frac{lk}{2}} =1$ is satisified if and only if $k$ is even. We observe that  for $l \geq 2$ 
 the equation 
\[ P(V,W) = q \cdot P(\lambda_1 V, \lambda_2 W)  \] 
can have nontrivial solutions only if the total degree of the polynomial is equal to $lk-2$ with even $k$. These solutions 
are of the form $(VW)^s$,  where $s = \frac{lk}{2} -1$. The expression $VW$ is a quadratic solution of the same equation.

To find $VW$ in terms of the variables $(v, w)$, all we need to do is to find all homogeneous quadratic solutions of the equation $p(v,w) =  q \cdot p( (q-1)v-w, v+w)$. 
Write $p(v,w) = c_1 v^2 + c_2 vw + c_3 w^2$.  Then we have the following system of equations for 
the coefficients $c_1, c_2, c_3$: 
\[ \left\{ \begin{array}{llll}  
(q(q-1)^2-1) c_1 & + q(q-1) c_2 &+ q c_3 &=0  \\ 
-2q(q-1) c_1  &+(q^2 -2q-1) c_2 &+ 2q c_3 &=0 \\ 
q c_1        &- q c_2         &+ (q-1)c_3   &= 0 
\end{array}  \right.  \]  

In particular, the system implies 
\[ \left\{ \begin{array}{ll} 
c_1(q^3 -q^2-1)  &= -c_3(q^2-q+1) \\ 
2c_1 (q^3 -q^2 -1) & = -c_2(q^2+1) 
\end{array} \right. \] 

If $(q^2-q+1) =0$, then $q$ is a primitive $6$th root of unity, and $c_1=0$, which immediately implies $c_2 = c_3=0$. If $(q^2+1)=0$, then $q$ is a primitive $4$th root of unity, and again $c_1=0$, and $c_2=c_3=0$. The expression 
$q^3-q^2-1 \neq 0$ for any root of unity $q$. Further,  
assuming $l \neq 4,6$, we can  express $c_3$ and $c_2$ in terms of $c_1$. Then the third equation    
leads to  the relation for $q$: 
\[  (q^2-1)(q^4 -q^3 +2q^2 +1)=  0 .\] 
We will consider three cases. 

\underline{$q=1$}.  \\ 
This is the case of the divergence that we considered before. The system of equations leads to the solution 
$c_1=c_2=c_3$, which results in a symmetric polynomial 
$p(v,w) = v^2  + vw + w^2$. The powers of $(v^2 + vw + w^2)$ are never antysymmetric and do not provide 
solutions for both conditions listed in the Proposition. Instead, the reason for the existence of nontrivial 
solutions in case of the divergence is that for $q=1$, we have 
$|\lambda_1| = |\lambda_2| =1$, and so the potential solutions of the polynomial 
equations are not restricted to the powers of $(VW)$. 

\underline{$q=-1$}. \\ 
This is the case of the superdivergence. In this case the system of equations for the coefficients yields the solution 
$c_1 = c_3$, $c_2 = 3 c_1$. 
Therefore, for each even degree $n$, the set of homogeneous polynomials of degree $n$ satisfying the condition 
$p(v,w) = -p(v+w, -v-2w)$ is empty if $n \equiv 0 \pmod 4$, and consists of the polynomials  proportional to 
$(v^2 + 3vw + w^2)^{\frac{n}{2}}$ if $n \equiv 2 \pmod 4$. 
All such polynomials are symmetric, and therefore adding the condition $p(v,w) =-p(w,v)$ implies that 
there are no homogeneous polynomial satisfies both conditions announced in the Proposition.  

\underline{$q \neq \pm 1$}. \\ 
This is the case of the $q$-divergence with $l \geq 3$. 
It is easy to see that solutions of the equation $q^4 -q^3 +2q^2 +1 =0$ are not roots of unity. Therefore, for $l \geq 3$, no 
quadratic homogeneous polynomial satisfies the equation $p(v,w) =  q \cdot p( (q-1)v-w, v+w)$, and consequently 
no homogeneous polynomial of any degree can satisfy this equation. 

Finally, we conclude that for $l \geq 2$, no homogeneous polynomial in two variables satisfies both conditions of the Proposition. 
\qed

Propositions \ref{qkern} and \ref{qker_poly} together provide a proof of Theorem \ref{qker}.

\section{Appendix} 

This section contains the computation of the matrix elements and characters of the representation $\rho_n$ 
of $D_{12} = \{s^2=1, r^6=1, srs^{-1}=r^{-1} \}$ on the space of degree $n$ homogeneous polynomials in two variables, 
defined by 

\[ g \cdot p(v,w) = p(g^{-1}(v,w)) ; \] 

\[  s \cdot  \begin{pmatrix} v \\ w \end{pmatrix} =  \begin{pmatrix} 0 &1 \\ 1 &0 \end{pmatrix}  \begin{pmatrix} v \\ w \end{pmatrix} = \begin{pmatrix}  w \\ v \end{pmatrix} \]
and
\[ r \cdot  \begin{pmatrix} v \\ w \end{pmatrix}  = \begin{pmatrix*}[r] 1 &1 \\ -1 & 0 \end{pmatrix*}   \begin{pmatrix} v \\ w \end{pmatrix} = \begin{pmatrix}  w + v \\ - v \end{pmatrix}  \]
 
\subsection{Matrix elements of the action}   
Let us compute the matrix elements $(\rho_n)_{ij}$ of this  representation, where the indices $i,j$ change from $1$ to $n+1$. 
 For the conjugacy classes of $e,s, r^3$ we have  

$\begin{array}{ll}

\rho_n(e)=  I_{n+1} & \chi_n(e) = n+1\vspace{0.3cm}\\ 

\rho_n(s) = \left( \begin{array}{ccc} 0& &1\\ \ &\ldots &\\ 1& &0 \end{array} \right)   & \chi_n(s) = \left\{ \begin{array}{l} 1 \;\; \text{if n even} \\ 0 \; \; \text{if n odd} \end{array} \right .\vspace{0.3cm}\\

\rho_n({r^3}) = \left\{\begin{array}{r} I_{n+1} \; \text{if n even} \\ -I_{n+1} \; \text{ if n odd} \end{array} \right. & \chi_n(r^3)= \left\{\begin{array}{r} n+1\; \text{ if n even} \\ -(n+1)\; \text{ if n odd} \end{array} \right.\vspace{0.3cm} \\
\end{array} $ 

Now let us compute the matrix coefficients of the action of $rs$. 
We will use the following notation for the binomial coefficients: 
\[  C_n^k = \frac{n!}{k! \cdot (n-k)!}  ,  \]
where $n \geq 0$ and $0 \leq k \leq n$. If $n$ and $k$ do not satisfy these restrictions, we will assume 
that $C_n^k =0$. 

Let $\{ a_0, a_1,  \ldots a_n \}$ denote the coefficients of a homogeneous polynomial of degree $n$, 
\[ p(x,y) = a_0 x^n + a_1 x^{n-1} y + \ldots + a_{n-1} x y^{n-1} + a_n y^n . \] 

The action of $rs$ is given by 
\[ \rho_n(rs) (p(x,y)) = p(x+y, -y)  = \sum_{j=0}^n a_j (x+y)^{n-j} (-y)^{j} = 
\sum_{j=0}^n a_j (-1)^j \sum_{i=0}^{n-j} C_{n-j}^i x^i y^{n-i} . \]
Changing the order of summation, we have 
\[ \rho_n(rs) (p(x,y)) = \sum_{i=0}^n  \left( \sum_{j=0}^{n-i} a_j (-1)^j  C_{n-j}^i \right) x^i y^{n-i} . \] 
The coefficients 
\[ (-1)^j C_{n-j}^i \] 
span one row of the matrix, the row that corresponds to $x^i y^{n-i}$. Since the first row of the matrix 
corresponds to $x^n$, and the $(n+1)$st row to $y^n$, the row corresponding to $x^i y^{n-i}$ has number $n-i+1$. 
Also, the index $j$ of the element in the row should change from $1$ to $n+1$ rather than from $0$ to $n$. Then we get: 
\[ ( \rho_n(rs) )_{i,j} = (-1)^{j-1} C_{n-j+1}^{n-i+1}   = (-1)^{j-1} C_{n-j+1}^{i-j} . \]
The last equality follows from the symmetry of the binomial coefficients. 

For the matrix elements of the action of $r$, we have 

\[ \rho_n(rs) \cdot \rho_n(s) = \rho_n(r) .\]  
Since $\rho_n(s)$ is the anti-diagonal matrix with $1$s on the anti-diagonal and zeros elsewhere, multiplication of $\rho_n(rs)$  by this matrix is equivalent to reading the rows of $\rho_n(rs)$ backwards. All we need to do to obtain $\rho_n(r)$ 
is to change variables, $j \to n+2-j$ in the expressions for the matrix elements of $\rho_n(rs)$: 
\[ ( \rho_n(r) )_{ij} = (\rho_n(rs))_{i, n+2-j} = (-1)^{n+1-j} C^{n+1-i}_{j-1} = (-1)^{n+1-j} C^{i+j -n-2}_{j-1} . \]  

For the action of $r^2$, we have: 
\[ \rho_n({r^2})(p(x,y)) = p(-x-y, x) = (-1)^n p(x+y, -x) .\]  
Therefore, the matrix of $r^2$ differs from the matrix of $rs$ by multiplication by $(-1)^n$ 
and reversing the order of elements in each column, which corresponds to the swap of variables $x \leftrightarrow y$. 
The reversing of the order of elements in each column is implemented by the change of variable: $ i \to n+2-i$. 
Then the matrix of $r^2$ is given by 
\[ (\rho_n(r^2))_{ij} = (-1)^n (\rho_n(rs))_{n+2-i, j} =  (-1)^{n+1-j}  C_{n-j+1}^{i-1} = (-1)^{n+1-j} C_{n-j+1}^{n+2-i-j} . \] 

\subsection{Characters of the representation $\rho_n$.} 

For $\chi_n(rs)$, we have 
\[ \chi_n(rs) = \sum_{i=1}^{n+1} (-1)^{i-1} C_{n-i-1}^{n-i-1} = \sum_{i=0}^n (-1)^n = \left[  
\begin{array}{ll} 1,   & \;\;\;  n \;\;{\rm even} \\  
                              0,  & \;\;\;  n \;\; {\rm odd}  \end{array} \right.   \] 

For $\chi_n(r)$, we compute using suitable changes of variables: 
\[ \chi_n(r) =  (-1)^n  \sum_{i=1}^{n+1}  (-1)^{i-1} C_{i-1}^{n-(i-1)}  \stackrel{k=i-1}{=} (-1)^n \sum_{k=0}^n (-1)^k C_k^{n-k} 
\stackrel{k \to n-k}{=} (-1)^n \sum_{k=0}^n (-1)^{n-k} C_{n-k}^k =  \] 
\[ = \sum_{k=0}^{\left[\frac{n}{2}\right] } (-1)^k C_{n-k}^k  . \]
The last equality holds because when $k > \left[ \frac{n}{2} \right]$ we have $k > n-k$ and $C_{n-k}^k =0$.

For $\chi_n(r^2)$, we compute 
\[ \chi_n(r^2) = (-1)^n \sum_{i=1}^{n+1} (-1)^{i-1} C_{n -(i-1)}^{i-1} \stackrel{k = i-1} = (-1)^n \sum_{k=0}^n (-1)^k C_{n-k}^k = 
(-1)^n \sum_{k=0}^{\left[\frac{n}{2}\right] } (-1)^k C_{n-k}^k  .\]

Therefore, we have 

$\begin{array}{ll} 
( \rho_n(rs) )_{ij} =  (-1)^{j-1} C_{n-j+1}^{n-i+1}  & \chi_n(rs) =  \sum_{i=1}^{n+1} (-1)^{i-1}= \left\{ \begin{array}{l} 1 \; \text{if n even} \\ 0 \; \text{if n odd} \end{array} \right.\vspace{0.3cm} \\
 
(\rho_n(r^2))_{ij} =  (-1)^{n+1-j}  C_{n-j+1}^{i-1} & \chi_n(r^2)  = (-1)^n \sum_{k=0}^{\left[\frac{n}{2}\right] } (-1)^k C_{n-k}^k \vspace{0.3cm}\\

( \rho_n(r) )_{ij}  = (-1)^{n+1-j} C^{n+1-i}_{j-1}  & \chi_n(r) = \sum_{k=0}^{\left[\frac{n}{2}\right] } (-1)^k C_{n-k}^k   \vspace{0.3cm}\\
\end{array} $

To compute  the  characters $\chi_n(r^2)$, $\chi_n(r)$ explicitly, we will use the following statement. 

\begin{lem} 
Let $n$ be a nonnegative integer, and denote 
\[ f_n = \sum_{k=0}^{[\frac{n}{2}]} (-1)^k C_{n-k}^k , \] 
where $C_{n-k}^k = \frac{(n-k)!}{k! (n-2k)!}$ is the binomial coefficient, and $\left[ \frac{n}{2} \right]$ denotes the integer part of 
$\frac{n}{2}$. 
Then the following recursive formula holds for all $n \geq 2$:  
\[ f_n = f_{n-1} - f_{n-2} . \]
In particular, the sequence $\{ f_n \}_{n \geq 0}$ is periodic with period $6$: 
\[ 1,1, 0, -1 -1, 0,  (1,1, 0, -1 -1, 0), \ldots . \]  
\end{lem} 

{\it Proof:} 

First suppose that $n \geq 2$ is even. Then we can write
\[ f_n = \sum_{k=0}^{\frac{n}{2}} (-1)^k C_{n-k}^k = 1 + \sum_{k=1}^{\frac{n}{2} -1} (-1)^k C_{n-k}^k + (-1)^{\frac{n}{2}} . \] 
Using the binomial coefficient identity  
\[ C_m^l  = C_{m-1}^{l-1} + C_{m-1}^l , \quad \quad 1 \leq l \leq m-1 ,\]  
we have 
\[ f_n =   1 + (-1)^{\frac{n}{2}} + \sum_{k=1}^{\frac{n}{2} -1} (-1)^k  C_{n-1-k}^{k} + \sum_{k=1}^{\frac{n}{2} -1} (-1)^k  C_{n-1-k}^{k-1}  = \] 
\[  = \left( 1 + \sum_{k=1}^{\frac{n}{2} -1} (-1)^k C_{(n-1)-k}^k  \right)  + \left( (-1)^{\frac{n}{2}} + (-1) \sum_{i=0}^{\frac{n}{2}-2} (-1)^i 
C_{(n-2)-i}^i  \right)  . \] 
The expression in the first parenthesis equals $\sum_{k=0}^{[\frac{n-1}{2}]} (-1)^k C_{n-1-k}^k = f_{n-1}$. In the second parenthesis, 
we have 
\[ (-1)^{\frac{n}{2}} =  (-1) \cdot (-1)^{\frac{n-2}{2}} = (-1)\cdot (-1)^{\frac{n-2}{2}}  
C_{\frac{n-2}{2}}^{\frac{n-2}{2}}, \] and this term can be included in the sum over $i$ for $i = \frac{n-2}{2}$: 
\[ f_n = \sum_{k=0}^{[\frac{n-1}{2}]} C_{(n-1)-k}^k  - \left( \sum_{i=0}^{\frac{n-2}{2}} (-1)^i C_{n-2-i}^i \right)  = f_{n-1} - f_{n-2} .\] 
The proof for an odd $n$ is similar.  Since $f_0 = f_1 =1$, we obtain the required sequence. 
\qed 

Then the values of $\chi_n$ on all conjugacy classes of $D_{12}$ are given in the table in the proof 
of Theorem \ref{dim_kern}.

\end{document}